# Subsets of Grassmannians Preserved by Mean Curvature Flows

Mu-Tao Wang

**Abstract**

Let $M = \Sigma_1 \times \Sigma_2$ be the product of two compact Riemannian manifolds of dimension $n \geq 2$ and two, respectively. Let $\Sigma$ be the graph of a smooth map $f : \Sigma_1 \mapsto \Sigma_2$, then $\Sigma$ is an $n$-dimensional submanifold of $M$. Let $\mathfrak{G}$ be the Grassmannian bundle over $M$ whose fiber at each point is the set of all $n$-dimensional subspaces of the tangent space of $M$. The Gauss map $\gamma : \Sigma \mapsto \mathfrak{G}$ assigns to each point $x \in \Sigma$ the tangent space of $\Sigma$ at $x$. This article considers the mean curvature flow of $\Sigma$ in $M$. When $\Sigma_1$ and $\Sigma_2$ are of the same non-negative curvature, we show a sub-bundle $\mathfrak{S}$ of the Grassmannian bundle is preserved along the flow, i.e. if the Gauss map of the initial submanifold $\Sigma$ lies in $\mathfrak{S}$, then the Gauss map of $\Sigma_t$ at any later time $t$ remains in $\mathfrak{S}$. We also show that under this initial condition, the mean curvature flow remains a graph, exists for all time and converges to the graph of a constant map at infinity . As an application, we show that if $f$ is any map from $S^n$ to $S^2$ and if at each point, the restriction of $df$ to any two dimensional subspace is area decreasing, then $f$ is homotopic to a constant map.

# 1 Introduction

The maximum principle has proved to be a powerful tool in partial differential equations. In particular, the maximum principle of parabolic systems for tensors developed by R. Hamilton [4] plays an important role in the study of geometric evolution equations. The guiding principle is the following: an invariant convex subset in the space of curvature tensors preserved by the associated ordinary differential equations is preserved by the parabolic partial differential equations. This has been applied to the study of Ricci flow and curvature flow of hypersurfaces. In this article, we apply this idea to higher codimension mean curvature flow.

Let $M = \Sigma_1 \times \Sigma_2$ be the product of two compact Riemannian manifolds of dimension $n \geq 2$ and two, respectively. Let $\Sigma$ be the graph of a map $f : \Sigma_1 \mapsto \Sigma_2$, then $\Sigma$ is an $n$-dimensional submanifold of $M$. Let $\mathfrak{G}$ be the Grassmannian bundle over $M$ whose fiber at each point is the set of all $n$-dimensional subspaces of the tangent space. The Gauss map $\gamma : \Sigma \mapsto \mathfrak{G}$ assigns to each point $x \in \Sigma$ the tangent space of $\Sigma$ at $x$. The tangent space of $M$ at $x$ splits as $T_{\pi_1(x)}\Sigma_1 \times T_{\pi_2(x)}\Sigma_2$. Let $\mathfrak{G}' \subset \mathfrak{G}$ be the sub-bundle consisted of the graphs of linear transformations from $T_{\pi_1(x)}\Sigma_1$ to $T_{\pi_2(x)}\Sigma_2$. We show there exists a sub-bundle $\mathfrak{S} \subset \mathfrak{G}'$ that is preserved along the mean curvature flow.

**Theorem A** *Let $M = \Sigma_1 \times \Sigma_2$ be the product of two compact flat Riemannian manifolds and suppose $\Sigma_2$ is two-dimensional. If the gauss map of a compact oriented submanifold $\Sigma$ of $M$ lies in $\mathfrak{S}$, then along the mean curvature flow the gauss map of $\Sigma_t$ remains in $\mathfrak{S}$. The flow of exists smoothly for all time and converges to a totally geodesic submanifold.*

This in particular implies $\Sigma_t$ is the graph of a map $f_t$. The sub-bundle $\mathfrak{S}$ is best described in terms of $f_t$. In fact, if we denote the singular values of $f_t$ by $\lambda_1$ and $\lambda_2$, then the Gauss map of $\Sigma_t$ lies in $\mathfrak{S}$ if and only if $|\lambda_1\lambda_2| < 1$.

When $\Sigma_1$ is of positive curvature, we prove the following.

**Theorem B** *Let $M = S^n(k_1) \times \Sigma_2$ be the product of a sphere of curvature $k_1 > 0$ and a two-dimensional compact Riemannian manifold $\Sigma_2$ of constant curvature $k_2$ and $k_1 \geq |k_2|$. If the gauss map of a compact oriented submanifold $\Sigma$ of $M$ lies in $\mathfrak{S}$, then along the mean curvature flow the gauss map of $\Sigma_t$ remains in $\mathfrak{S}$. The flow exists smoothly for all time and converges to a totally geodesic submanifold.*

Theorem A and B are proved by calculating the evolution equations of the Gauss map and applying maximum principle. The prototype is the following equation in the hypersurface case

$$\frac{d}{dt}N = \Delta N + |A|^2 N$$

where $N$ denotes the unit normal vector and $|A|^2$ is the norm of the second fundamental form. If we take inner product of $N$ with a constant vector $\nu$, it is not hard to see that $\min_{\Sigma_t} < N, \nu >$ is non-decreasing in time. This is one of the key observation in [1] and [2] where the mean curvature flow of of entire graph of codimension one was studied. In codimension one case, $N$ contains all the information of the Gauss map. While in higher codimension, a whole parabolic systems is needed in order to describe the evolution of the Gauss map.

The following is an application to higher homotopy groups of $S^2$.

**Corollary** *If $f$ is any map from $S^n$ to $S^2$, $n \geq 2$ and if at each point, the restriction of $df$ to any two dimensional subspace is area decreasing, then $f$ is homotopic to a constant map along the mean curvature flow.*

When $n = 2$, this is the same as saying the Jacobian of $f$ is less than 1. In this case, $f$ is of degree 0 and thus homotopic to a constant map. This homotopy can be realized through the mean curvature flow as was proved in [8]. As a contrast, the standard Hopf map from $S^3$ to $S^2$ has $|\lambda_1 \lambda_2| = 4$.

I am indebted to Professor D. H. Phong and Professor S.-T. Yau for their constant encouragement and unending support. I have benefitted greatly from the conversation I have with Professor R. Hamilton and Professor M-P Tsui.

# 2 Analysis of Grassmannian bundle

Let us first describe the sub-bundle $\mathfrak{S}$. Let $V_1$ be an $n$-dimensional inner product space and $V_2$ a two-dimensional inner product space. Let $G(n, n+2)$ be the Grassmannian of all $n$-dimensional subspaces of $V_1 \times V_2$. Let $G' \subset G(n, n+2)$ be the set of all $n$-dimensional subspaces that can be written as graphs over $V_1$. For any $P \in G'$, $P$ is the graph of a linear transformation $\mathfrak{p} : V_1 \mapsto V_2$. Then $(\mathfrak{p})^T \mathfrak{p}$ is self-adjoint and thus diagonalizable. The eigenvalues are denoted by $\{\lambda_1^2, \lambda_2^2\}$. $\lambda_1$ and $\lambda_2$ are the singular values of $\mathfrak{p}$. We now define $S$.

$$S = \{P \in G' |\, 1 - |\lambda_1 \lambda_2| > 0\}$$

This is equivalent to saying $\mathfrak{p}$ is area decreasing on any two dimensional subspace of $V_1$.

Now let $M$ be the product of two Riemannian manifolds $\Sigma_1 \times \Sigma_2$ of dimension $n$ and 2 respectively. Let $\mathfrak{G}$ be the Grassmannian bundle on $M$ whose fibers are isomorphic to $G(n, n+2)$. At each point $x$, $T_xM$ splits as the product of $T_{\pi_1(x)}\Sigma_1$ and $T_{\pi_2(x)}\Sigma_2$. The Riemannian structures on $T_{\pi_1(x)}\Sigma_1 = V_1$ and $T_{\pi_2(x)}\Sigma_2 = V_2$ defines the subset $S$ of the fiber of $\mathfrak{G}$ at $x$.

**Definition 2.1** *$\mathfrak{S}$ is the sub-bundle of the Grassmannian bundle $\mathfrak{G}$ whose fiber at each point consists of $S$.*

Let $\Sigma$ be the graph of a smooth map $f : \Sigma_1 \mapsto \Sigma_2$. $T_x\Sigma$ is the graph of the differential of $f$ at $x$, $df : T_{\pi_1(x)}\Sigma_1 \mapsto T_{\pi_2(x)}\Sigma_2$. Notice that we abuse the notation so that $T_x\Sigma, T_{\pi_1(x)}\Sigma_1$ and $T_{\pi_2(x)}\Sigma_2$ all denote subspaces of $T_xM$. At any point $x$, let $\lambda_1$,$\lambda_2$ be the singular values of $df$. They are well-defined up to a sign. Define

$$\eta(x) = \frac{1 - |\lambda_1 \lambda_2|}{\sqrt{(1+\lambda_1^2)(1+\lambda_2^2)}}$$

$\eta$ is a function on $\Sigma$.

**Proposition 2.1** *$\eta > 0$ on $\Sigma$ if and only if the Gauss map of $\Sigma$ lies in $\mathfrak{S}$.*

Later we shall give a characterization of $\eta$ in terms of differential forms on $M$. Any differential form $\Omega$ on a Riemannian manifold can be considered as a function on the Grassmannian bundle $\mathfrak{G}$ of appropriate dimension. The *comass* of $\Omega$ at $x$ is defined to be the supremum of $\Omega$ on $\mathfrak{G}_x$, the fiber of $\mathfrak{G}$ at $x$. This is an important concept in calibrated geometry, see Federer [3] or Harvey-Lawson [5]. Another description of $\mathfrak{S}$ can be given in terms of the comass.

**Proposition 2.2** *If $\Sigma$ is the graph of $f : \Sigma_1 \mapsto \Sigma_2$, then the Gauss map of $\Sigma$ lies in $\mathfrak{S}$ if and only if the comass of $f^*\Omega_2$ is less than one.*

Here $\Omega_2$ is the volume form on $\Sigma_2$ and $f^*\Omega_2$ is considered as a 2-form on $\Sigma_1$. Of course the comass is taken over all two-dimensional subspaces of the tangent space of $\Sigma_1$.

## 3 Evolution equation of $n$ form

In this section, we calculate the evolution equation of the restriction of an $n$-form to an $n$-dimensional submanifold moving by the mean curvature flow. The case for a parallel form was calculated in [10]. Here we need to keep track of the terms that involve covariant derivatives of $\Omega$.

We assume $M$ is an $n+m$ dimensional Riemannian manifold with an $n$ form $\Omega$. Let $F : \Sigma \mapsto M$ be an isometric immersion of an $n$-dimensional submanifold. We shall compute near a point $p \in \Sigma$. We choose arbitrary orthonormal frames $\{e_i\}_{i=1\cdots n}$ for $T\Sigma$ and $\{e_\alpha\}_{\alpha=n+1,\cdots,n+m}$ for $N\Sigma$. $\nabla^M$ denotes the covariant derivative on $M$ and $\nabla^\Sigma$ denotes the covariant derivative on $\Sigma$, which is simply the tangent part of $\nabla^M$. $\nabla^M\Omega$ is the covariant derivative of $\Omega$ on $M$ and $\nabla^\Sigma\Omega$ will denote the covariant derivative of the restriction of $\Omega$ to $\Sigma$.

We first calculate the covariant derivative of the restriction of $\Omega$ on $\Sigma$.

$$
\begin{aligned}
&(\nabla^\Sigma_{e_k}\Omega)(e_{i_1},\cdots,e_{i_n})\\
&= e_k(\Omega(e_{i_1},\cdots,e_{i_n})) - \Omega(\nabla^\Sigma_{e_k}e_{i_1},\cdots,e_{i_n}) - \cdots - \Omega(e_{i_1},\cdots,\nabla^\Sigma_{e_k}e_{i_n})\\
&= (\nabla^M_{e_k}\Omega)(e_{i_1},\cdots,e_{i_n}) + \Omega(\nabla^M_{e_k}e_{i_1} - \nabla^\Sigma_{e_k}e_{i_1},\cdots,e_{i_n}) + \cdots + \Omega(e_{i_1},\cdots,\nabla^M_{e_k}e_{i_n} - \nabla^\Sigma_{e_k}e_{i_n})
\end{aligned}
$$

This equation can be abbreviated using the second fundamental form of $F$, $h_{\alpha ij} =< \nabla^M_{e_i}e_j, e_\alpha >$.

$$
\Omega_{i_1\cdots i_n,k} = (\nabla^M_{e_k}\Omega)(e_{i_1},\cdots,e_{i_n}) + \Omega_{\alpha i_2\cdots i_n}h_{\alpha i_1 k} + \cdots + \Omega_{i_1\cdots i_{n-1}\alpha}h_{\alpha i_n k} \quad (3.1)
$$

Likewise, in $\Omega(e_\alpha, e_{i_2},\cdots,e_{i_n})$, $\Omega$ is considered as a section of $(N\Sigma)^* \wedge (\wedge(T\Sigma)^*)$.

$$
\Omega_{\alpha i_2\cdots i_n,k} = (\nabla^M_{e_k}\Omega)(e_\alpha, e_{i_1},\cdots,e_{i_n}) - \Omega_{j i_2\cdots i_n}h_{\alpha jk} + \Omega_{\alpha\beta i_3\cdots i_n}h_{\beta i_2 k} + \cdots + \Omega_{\alpha i_2\cdots i_{n-1}\beta}h_{\beta i_n k} \quad (3.2)
$$

Now we calculate the second covariant derivative of the restriction of $\Omega$ on $\Sigma$.

$$
\begin{aligned}
&(\nabla^\Sigma_{e_k}\nabla^\Sigma_{e_k}\Omega)(e_1,\cdots e_n)\\
&= e_k((\nabla^\Sigma_{e_k}\Omega)(e_1,\cdots e_n)) - (\nabla^\Sigma_{e_k}\Omega)(\nabla^\Sigma_{e_k}e_1,\cdots,e_n) - \cdots - (\nabla^\Sigma_{e_k}\Omega)(e_1,\cdots,\nabla^\Sigma_{e_k}e_n)
\end{aligned}
$$

The term $(\nabla^{\Sigma}_{e_k}\Omega)(\nabla^{\Sigma}_{e_k}e_1,\cdots,e_n)$ equals zero because $\nabla^{\Sigma}_{e_k}e_1$ is a tangent vector perpendicular to $e_1$ and thus a linear combination of $e_2,\cdots,e_n$. Likewise, other similar terms vanish.

$$
\begin{aligned}
&(\nabla^{\Sigma}_{e_k}\nabla^{\Sigma}_{e_k}\Omega)(e_1,\cdots e_n)\\
=&e_k[(\nabla^{M}_{e_k}\Omega)(e_1,\cdots,e_n)+\Omega_{\alpha 2\cdots n}h_{\alpha 1k}+\cdots+\Omega_{1\cdots n-1\alpha}h_{\alpha nk}]\\
=&(\nabla^{M}_{e_k}\nabla^{M}_{e_k}\Omega)(e_1,\cdots,e_n)+(\nabla^{M}_{e_k}\Omega)(\nabla^{M}_{e_k}e_1,\cdots,e_n)+\cdots+(\nabla^{M}_{e_k}\Omega)(e_1,\cdots,\nabla^{M}_{e_k}e_n)\\
&+\Omega_{\alpha 2\cdots n,k}h_{\alpha 1k}+\cdots+\Omega_{1\cdots n-1\alpha,k}h_{\alpha nk}\\
&+\Omega_{\alpha 2\cdots n}h_{\alpha 1k,k}+\cdots+\Omega_{1\cdots n-1\alpha}h_{\alpha nk,k}
\end{aligned}
\tag{3.3}
$$

Now $\nabla^{M}_{e_k}e_i=h_{\alpha ik}e_\alpha+\nabla^{\Sigma}_{e_k}e_i$ and $(\nabla^{M}_{e_k}\Omega)(\nabla^{\Sigma}_{e_k}e_1,\cdots,e_n)=0$
Therefore,

$$
\begin{aligned}
\Omega_{1\cdots n,kk}=&(\nabla^{M}_{e_k}\nabla^{M}_{e_k}\Omega)(e_1,\cdots,e_n)\\
&+(\nabla^{M}_{e_k}\Omega)(e_\alpha,\cdots,e_n)h_{\alpha 1k}+\cdots+(\nabla^{M}_{e_k}\Omega)(e_1,\cdots,e_\alpha)h_{\alpha nk}\\
&+\Omega_{\alpha 2\cdots n,k}h_{\alpha 1k}+\cdots+\Omega_{1\cdots n-1\alpha,k}h_{\alpha nk}\\
&+\Omega_{\alpha 2\cdots n}h_{\alpha 1k,k}+\cdots+\Omega_{1\cdots n-1\alpha}h_{\alpha nk,k}
\end{aligned}
\tag{3.4}
$$

Plug equation (3.2) into (3.4) and apply the Codazzi equation $h_{\alpha ki,k}=h_{\alpha,i}+R_{\alpha kki}$ where $R$ is the curvature operator of $M$.

$$
\begin{aligned}
(\Delta^{\Sigma}\Omega)_{1\cdots n}=&-\Omega_{12\cdots n}\sum_{\alpha,k}(h^2_{\alpha 1k}+\cdots+h^2_{\alpha nk})\\
&+2\sum_{\alpha,\beta,k}[\Omega_{\alpha\beta 3\cdots n}h_{\alpha 1k}h_{\beta 2k}+\Omega_{\alpha 2\beta\cdots n}h_{\alpha 1k}h_{\beta 3k}+\cdots+\Omega_{1\cdots(n-2)\alpha\beta}h_{\alpha(n-1)k}h_{\beta nk}]\\
&+\sum_{\alpha,k}\Omega_{\alpha 2\cdots n}h_{\alpha,1}+\cdots+\Omega_{1\cdots(n-1)\alpha}h_{\alpha,n}\\
&+\sum_{\alpha,k}\Omega_{\alpha 2\cdots n}R_{\alpha kk1}+\cdots+\Omega_{1\cdots(n-1)\alpha}R_{\alpha kkn}\\
&+(\nabla^{M}_{e_k}\nabla^{M}_{e_k}\Omega)(e_1,\cdots,e_n)\\
&+2(\nabla^{M}_{e_k}\Omega)(e_\alpha,\cdots,e_n)h_{\alpha 1k}+\cdots+2(\nabla^{M}_{e_k}\Omega)(e_1,\cdots,e_\alpha)h_{\alpha nk}
\end{aligned}
\tag{3.5}
$$

We notice that $(\Delta^\Sigma \Omega)_{1\cdots n} = \Delta(\Omega(e_1, \cdots, e_n))$, where the $\Delta$ on the right hand side is the Laplacian of functions on $\Sigma$.

The terms in the bracket are formed in the following way. Choose two different indexes from 1 to $n$, replace the smaller one by $\alpha$ and the larger one by $\beta$. There are a total of $\frac{n(n-1)}{2}$ such terms.

Now we consider when $\Sigma = \Sigma_t$ is a time slice of a mean curvature flow in $M$ by $\frac{d}{dt}F_t = H_t$. Notice that here we require the velocity vector is in the normal direction. We can extend $e_1, \cdots, e_n$ to a local coordinate $\{\partial_i = \frac{\partial}{\partial x^i}\}$ on $\Sigma$, then

$$
\begin{aligned}
&\frac{d}{dt}\Omega(\partial_1, \cdots, \partial_n) \\
&= (\nabla^M_H \Omega)(\partial_1, \cdots, \partial_n) + \Omega((\nabla_{\partial_1} H)^N, \partial_2, \cdots, \partial_n) + \cdots + \Omega(\partial_1, \partial_2, \cdots, (\nabla_{\partial_n} H)^N) \\
&+ \Omega((\nabla_{\partial_1} H)^T, \partial_2, \cdots, \partial_n) + \cdots + \Omega(\partial_1, \partial_2, \cdots, (\nabla_{\partial_n} H)^T)
\end{aligned}
$$

Since $\frac{d}{dt}g_{ij} =< (\nabla_{\partial_i} H)^T, \partial_j >$, if we choose a orthonormal frame and evolve the frame with respect to time so that it remains orthonormal, the terms in the last line vanish.

$$
\frac{d}{dt} * \Omega = *(\nabla^M_H \Omega) + \Omega_{\alpha 2\cdots n} h_{\alpha,1} + \cdots + \Omega_{1\cdots(n-1)\alpha} h_{\alpha,n}
$$

Combine this with equation (3.5) we get the parabolic equation satisfied by $*\Omega$.

**Proposition 3.1** *If $\Sigma_t$ is a time slice of an $n$-dimensional mean curvature flow in $M$ and $\Omega$ is an $n$-form on $M$. For any point $p \in \Sigma_t$, let $\{e_1, \cdots, e_n\}$ be an orthonormal frame of $T\Sigma_t$ near $p$ and $\{e_{n+1}, \cdots, e_{n+m}\}$ be an orthonormal* *frame of the normal bundle of $\Sigma_t$ near $p$. Then $*\Omega = \Omega(e_1, \cdots, e_n)$ satisfies*

$$\begin{aligned}\frac{d}{dt} * \Omega = \Delta * \Omega + *\Omega(\sum_{\alpha,i,k} h^2_{\alpha ik})\\ - 2\sum_{\alpha,\beta,k}[\Omega_{\alpha\beta 3\cdots n} h_{\alpha 1k} h_{\beta 2k} + \Omega_{\alpha 2\beta\cdots n} h_{\alpha 1k} h_{\beta 3k} + \cdots + \Omega_{1\cdots(n-2)\alpha\beta} h_{\alpha(n-1)k} h_{\beta nk}]\\ - \sum_{\alpha,k}[\Omega_{\alpha 2\cdots n} R_{\alpha kk1} + \cdots + \Omega_{1\cdots(n-1)\alpha} R_{\alpha kkn}]\\ + *(\nabla^M_H \Omega) - (\nabla^M_{e_k}\nabla^M_{e_k}\Omega)(e_1, \cdots, e_n)\\ - 2(\nabla^M_{e_k}\Omega)(e_\alpha, \cdots, e_n) h_{\alpha 1k} - \cdots - 2(\nabla^M_{e_k}\Omega)(e_1, \cdots, e_\alpha) h_{\alpha nk}\end{aligned} \tag{3.6}$$

*where $\Delta$ denotes the time-dependent Laplacian on $\Sigma_t$.*

# 4 Proof of Theorem

Let us prove Theorem A now. We recall the statement.

**Theorem A** *Let $M = \Sigma_1 \times \Sigma_2$ be the product of two compact flat Riemannian manifolds of dimension $n$ and 2 respectively. If the gauss map of a compact oriented submanifold $\Sigma$ of $M$ lies in $\mathfrak{S}$, then along the mean curvature flow the gauss map of $\Sigma_t$ remains in $\mathfrak{S}$. The flow exists smoothly for all time and converges to a totally geodesic submanifold.*

*Proof.* Let $\Sigma_t$ be the mean curvature flow of $\Sigma$ given by a family of immersions $F : \Sigma \times [0, T) \mapsto M$. In the following calculation, it is useful to consider the total space of the mean curvature flow as $\Sigma \times [0, T)$. At each instant $t$, $\Sigma$ is equipped with the induced metric by $F_t$. All geometric quantities defined on the image of $F(\cdot, t)$ are considered as defined on $\Sigma$.

Let $\Omega_1$ and $\Omega_2$ be the volume form of $\Sigma_1$ and $\Sigma_2$ respectively. They can be considered as parallel forms on $M$. Suppose initially the image of the Gauss map of $\Sigma$ is in $\mathfrak{S}$. We may assume $\Sigma$ is the graph of a map $f : \Sigma_1 \mapsto \Sigma_2$. This implies $\eta_1 = *\Omega_1 > 0$ and $\eta > 0$ on $\Sigma$ at $t = 0$ by Proposition 2.1.

We shall characterize $\eta$ in terms of differential forms. Consider $\Xi$ the collection of $n$ forms on $M$ of the following type.

$$\Xi = \{\, \Omega_1 - \Omega_2 \wedge \omega \,|\, \omega \text{ is any parallel simple } (n-2) \text{ form of comass one on } \Sigma_1\}$$

At any point $x$, by Singular Value Decomposition we can take an orthonormal basis $\{a_i\}_{i=1\cdots n}$ for $T_{\pi_1(x)}\Sigma_1$ and $\{a_\alpha\}_{\alpha=n+1,n+2}$ for $T_{\pi_2(x)}\Sigma_2$ so that $df(a_i) = \lambda_i a_{n+i}$, $a_1^* \wedge \cdots \wedge a_n^*$ is the volume form of $T_{\pi_1(x)}\Sigma_1$ and $a_{n+1}^* \wedge a_{n+2}^*$ is the volume form for $T_{\pi_2(x)}\Sigma_2$ . Therefore,

$$\{e_1 = \frac{1}{\sqrt{1+\lambda_1^2}}(a_1 + \lambda_1 a_{n+1}), e_2 = \frac{1}{\sqrt{1+\lambda_2^2}}(a_2 + \lambda_2 a_{n+2}), e_3 = a_3, \cdots, e_n = a_n\} \tag{4.1}$$

forms an orthonormal basis for $T_x\Sigma$ and

$$\{e_{n+1} = \frac{1}{\sqrt{1+\lambda_1^2}}(a_{n+1} - \lambda_1 a_1), e_{n+2} = \frac{1}{\sqrt{1+\lambda_2^2}}(a_{n+2} - \lambda_2 a_2)\} \tag{4.2}$$

an orthonormal basis for $N_x\Sigma$. Thus,

$$\begin{aligned} *(\Omega_1 - \Omega_2 \wedge \omega) &= (\Omega_1 - \Omega_2 \wedge \omega)(e_1, \cdots, e_n) \\ &= \frac{1}{\sqrt{(1+\lambda_1^2)(1+\lambda_2^2)}}(1 - \lambda_1\lambda_2\omega(a_3, \cdots, a_n)) \end{aligned} \tag{4.3}$$

On the other hand

$$\eta_1 = \frac{1}{\sqrt{(1+\lambda_1^2)(1+\lambda_2^2)}}$$

Recall

$$\eta = \frac{1 - |\lambda_1\lambda_2|}{\sqrt{(1+\lambda_1^2)(1+\lambda_2^2)}}$$

It is not hard to see

$$\eta(x) = \min_{\Omega\in\Xi} *\Omega(x)$$

Suppose at $t = t_0$, the image of the Gauss map hits the boundary of $\mathfrak{G}$ for the first time. Therefore each $\Sigma_t$, $t < t_0$ can be written as the graph of $f_t : \Sigma_1 \mapsto \Sigma_2$ and the singular values of $f_t$ satisfy $|\lambda_1\lambda_2| < 1$.

We claim $\Sigma_{t_0}$ remains a graph. Indeed, since $\Omega_1$ is a parallel form, $\eta_1$ satisfies the following equation by equation ( 3.6).

$$\frac{d}{dt}\eta_1 = \Delta\eta_1 + \eta_1[\sum_{\alpha,i,k} h^2_{\alpha ik} - 2\sum_k \lambda_1\lambda_2(h_{n+1,1k}h_{n+2,2k} - h_{n+2,1k}h_{n+1,2k})] \tag{4.4}$$

where we use

$$\Omega_1(e_{n+1}, e_{n+2}, e_3, \cdots, e_n) = \frac{\Omega_1(a_{n+1} - \lambda_1 a_1, a_{n+2} - \lambda_2 a_2, a_3, \cdots, a_n)}{\sqrt{(1+\lambda_1^2)(1+\lambda_2^2)}} = \lambda_1\lambda_2\eta_1$$

Notice this equation is valid at any point $x$. Since $|\lambda_1\lambda_2| < 1$ for $0 \leq t < t_0$, applying maximum principle to equation (4.4) implies $\min_{\Sigma_t} \eta_1$ is non-decreasing in $t$ and thus $\eta_1 > 0$ at $t_0$.

Now $\eta$ is well-defined at $t_0$. Take any $p$ so that $\eta(p, t_0) = 0$, we shall show that $\frac{d}{dt}|_{t=t_0}\eta \geq 0$ at $p$.

It is clear that $\lambda_1\lambda_2 \neq 0$ at $p$. Otherwise, $\eta_1 = \eta = 0$, a contradiction.

By the previous characterization of $\eta$ and Hamilton's maximum principle [4] , we only need to show $\frac{d}{dt}|_{t=t_0} * \Omega \geq 0$ at the point $p$ for any $\Omega \in \Xi$ such that $*\Omega(p) = \eta(p)$. At $p$, we apply Singular Value Decomposition to get an orthonormal basis $\{a_i\}_{i=1,\cdots n}$ for $T_{\pi_1(p)}\Sigma_1$ as before. Such $\Omega$ is of the form $\Omega_1 - \Omega_2 \wedge \omega$ with $\omega(a_3, \cdots, a_n) = 1$ or $\omega = a_3^* \wedge \cdots \wedge a_n^*$ by equation (4.3).

$*\Omega$ satisfies

$$\begin{aligned}\frac{d}{dt} * \Omega = {} & \Delta * \Omega + *\Omega(\sum_{\alpha,i,k} h^2_{\alpha ik}) \\ & - 2\sum_{\alpha,\beta,k}[\Omega_{\alpha\beta 3\cdots n}h_{\alpha 1k}h_{\beta 2k} + \Omega_{\alpha 2\beta\cdots n}h_{\alpha 1k}h_{\beta 3k} + \cdots + \Omega_{1\cdots(n-2)\alpha\beta}h_{\alpha(n-1)k}h_{\beta nk}]\end{aligned} \tag{4.5}$$

At this point $p$,

$$(\Omega_1 - \Omega_2 \wedge \omega)(e_{n+1}, e_{n+2}, e_{i_1}, \cdots e_{i_{n-2}}) = \frac{1}{\sqrt{(1+\lambda_1^2)(1+\lambda_2^2)}}(\lambda_1\lambda_2 - 1)\omega(e_{i_1}, \cdots, e_{i_{n-2}})$$

Thus

$$\frac{d}{dt} * \Omega = \Delta * \Omega + *\Omega[|A|^2 + 2(h_{n+1,1k}h_{n+2,2k} - h_{n+1,2k}h_{n+2,1k})]$$

This can be completed square and we get

$$\frac{d}{dt} * \Omega = \Delta * \Omega + *\Omega[\sum_{\alpha,2<i\leq n,k} h_{\alpha ik}^2 + \sum_k (h_{n+1,1k} + h_{n+2,2k})^2 + \sum_k (h_{n+1,2k} - h_{n+2,1k})^2]$$

Therefore $\frac{d}{dt}*\Omega \geq 0$ at $(p, t_0)$. Since this is true for any $\Omega$ that achieves the minimum of $*\Omega$ in $\Xi$, we have $\frac{d}{dt}\eta \geq 0$. Thus the sub-bundle $\mathfrak{S}$ is preserved along the mean curvature flow.

Now we prove long time existence and convergence. By a similar argument, we can show if $\min \eta = \delta > 0$ at $t = 0$, then this is preserved along the flow. This implies in particular,

$$|\lambda_1 \lambda_2| \leq 1 - \delta, \tag{4.6}$$

and

$$\sqrt{(1 + \lambda_1^2)(1 + \lambda_2^2)} \leq \frac{1}{\delta} \tag{4.7}$$

Since $|\lambda_1 \lambda_2| \leq 1 - \delta$, we have

$$\frac{d}{dt} * \Omega_1 \geq \Delta * \Omega_1 + \delta * \Omega_1 |A|^2 \tag{4.8}$$

In particular, $*\Omega_1$ has a uniform lower bound, each $\Sigma_t$ can be written as the graph of a map $f_t : \Sigma_1 \mapsto \Sigma_2$, and $f_t$ has uniform gradient bound.

Integrating $\frac{d}{dt} * \Omega_1 \geq \Delta * \Omega_1 + \delta * \Omega_1 |A|^2$ over space and time from $t = 0$ to $t = \infty$ we get

$$\delta \int_0^\infty \int_{\Sigma_t} *\Omega_1 |A|^2 \leq \int_0^\infty \int_{\Sigma_t} *\Omega_1 |H|^2$$

For a mean curvature flow, $\int_0^\infty \int_{\Sigma_t} |H|^2 < \infty$, thus $\int_0^\infty \int_{\Sigma_t} |A|^2 < \infty$. We can extract a subsequence $t_i \to \infty$ such that $\int_{\Sigma_{t_i}} |A|^2 \to 0$. Because each $f_i$ has bounded gradient, this is the same as $\int_{\Sigma_1} |\nabla df_i|^2 \to 0$. Therefore $df_i$

is in $W^{1,2}$ which is compactly imbedded in $L^{\frac{2n}{n-2}}$. We can further extract a convergent subsequence which converges in $C^\alpha \cap W^{1,\frac{2n}{n-2}}$ norm. Apply the Sobolev inequality shows $df_i$ converges to a constant and the limit of $f_i$ is a totally geodesic submanifold. The uniform convergence of $f_t$ follows as the proof of Theorem C in [9], which uses the property that distance function to any totally geodesic submanifold in a Riemannian manifold of non-positive sectional curvature is convex. We remark that in this case, the limit is actually the graph of a linear map.

□

**Theorem B** *Let $M = S^n(k_1) \times \Sigma_2^m$ be the product of a sphere of curvature $k_1 > 0$ and a compact Riemannian manifold $\Sigma_2$ of constant curvature $k_2$ and $k_1 \geq |k_2|$. If the gauss map of a compact oriented submanifold $\Sigma^n$ of $M$ lies in $\mathfrak{S}$, then along the mean curvature flow the gauss map of $\Sigma_t^n$ remains in $\mathfrak{S}$. The flow exists smoothly for all time and converges to an $S^n(k_1)$ slice in $M$.*

*Proof.*

The proof follows the same strategy as that of Theorem A. We actually show the image of the Gauss map never hits the boundary of $\mathfrak{S}$. Suppose the contrary happens at $t = t_0$. Again, we look at the equation of $*\Omega_1$. Using $|\lambda_1\lambda_2| \leq 1$ for $0 \leq t < t_0$ and Proposition 3.2 in [10], we see

$$\begin{aligned} &\frac{d}{dt} * \Omega_1 \geq \Delta * \Omega_1 \\ &+ *\Omega_1 \sum_i \frac{\lambda_i^2}{1+\lambda_i^2}[k_1(\sum_{j\neq i}\frac{2}{1+\lambda_j^2}) + k_2(1-n+\sum_{j\neq i}\frac{2}{1+\lambda_j^2})] \end{aligned} \tag{4.9}$$

The last term comes from the curvature of $M$. Rewrite

$$k_1(\sum_{j\neq i}\frac{2}{1+\lambda_j^2})+k_2(1-n+\sum_{j\neq i}\frac{2}{1+\lambda_j^2}) = \frac{k_1-k_2}{2}(n-1)+\frac{k_1+k_2}{2}(\sum_{j\neq i}\frac{2}{1+\lambda_j^2}+1-n)$$

We claim the curvature term is always non-negative under our assumption. Since $k_1 - k_2 \geq 0, k_1 + k_2 \geq 0$, we only need to show

$$\sum_{i=1}^n \frac{\lambda_i^2}{1+\lambda_i^2}(1-n+\sum_{k\neq i}^n \frac{2}{1+\lambda_k^2}) \geq 0$$

This is indeed

$$\frac{\lambda_1^2}{1+\lambda_1^2}(n-3+\frac{2}{1+\lambda_2^2})+\frac{\lambda_2^2}{1+\lambda_2^2}(n-3+\frac{2}{1+\lambda_1^2})$$

This can be rewritten as

$$(n-2)\frac{\lambda_1^2+\lambda_2^2+2\lambda_1^2\lambda_2^2}{(1+\lambda_1^2)(1+\lambda_2^2)}+\frac{\lambda_1^2+\lambda_2^2-2\lambda_1^2\lambda_2^2}{(1+\lambda_1^2)(1+\lambda_2^2)}$$

which is non-negative under the assumption $|\lambda_1\lambda_2| \leq 1$.

Therefore, at $t_0$, $\Sigma_{t_0}$ remains a graph and $\min_{\Sigma_{t_0}} \eta = 0$. Take any $p$ so that $\eta(p, t_0) = 0$. We may assume $\lambda_1 > 0, \lambda_2 > 0$ at $p$.

As before, we choose orthonormal basis at $p$ that corresponds to the singular value decomposition of $df$. We can extend the orthonormal basis $\{a_i\}_{i=1\cdots n}$ at $T_{\pi_1(p)}\Sigma_1$ to an orthonormal frame field in a neighborhood $U_1 \subset \Sigma_1$ such that at this point $p$, $\nabla^M a_i = 0$, $i = 1 \cdots n$. This is possible because the Riemannian structure is a product and each $\Sigma_1$ slice is totally geodesic in $M$. Take

$$\Omega = \Omega_1 - \Omega_2 \wedge \omega$$

where $\omega = a_3^* \wedge \cdots \wedge a_n^*$.

$\Omega$ is an $n$ form defined on $U = U_1 \times \Sigma_2$ that satisfies $\nabla^M \Omega = 0$ at $p$.

Now we extend $\Omega$ to a global form on $M$. Take a cut-off function $\phi$ such that $\phi \equiv 1$ in a neighborhood of $p$ and $\phi$ has compact support. Then

$$\Omega = \Omega_1 - \phi\Omega_2 \wedge \omega$$

is such a global extension. Now $*\Omega(p, t_0) = 0$ and $(p, t_0)$, $*\Omega > 0$ for $0 \leq t < t_0$ and $*\Omega \geq 0$ at $t_0$. Therefore $\frac{d}{dt} * \Omega \leq 0$ and $\Delta * \Omega \geq 0$ at$(p, t_0)$.

We recall the evolution equation from equation (3.6).

$$\begin{aligned}
\frac{d}{dt} * \Omega = {} & \Delta * \Omega + *\Omega(\sum_{\alpha,i,k} h^2_{\alpha ik}) \\
& - 2\sum_{\alpha,\beta,k}[\Omega_{\alpha\beta 3\cdots n} h_{\alpha 1k} h_{\beta 2k} + \Omega_{\alpha 2\beta\cdots n} h_{\alpha 1k} h_{\beta 3k} + \cdots + \Omega_{1\cdots(n-2)\alpha\beta} h_{\alpha(n-1)k} h_{\beta nk}] \\
& - \sum_{\alpha,k}[\Omega_{\alpha 2\cdots n} R_{\alpha kk1} + \cdots + \Omega_{1\cdots(n-1)\alpha} R_{\alpha kkn}] \\
& + *(\nabla^M_H \Omega) - (\nabla^M_{e_k}\nabla^M_{e_k}\Omega)(e_1, \cdots, e_n) \\
& - 2(\nabla^M_{e_k}\Omega)(e_\alpha, \cdots, e_n) h_{\alpha i_1 k} - \cdots - 2(\nabla^M_{e_k}\Omega)(e_1, \cdots, e_\alpha) h_{\alpha i_n k}
\end{aligned} \tag{4.10}$$

By the way $\Omega$ is constructed, $\nabla^M\Omega = 0$ at $p$. We claim the term $(\nabla^M_{e_k}\nabla^M_{e_k}\Omega)(e_1, \cdots, e_n)$ also vanishes at $p$.

In fact, consider

$$\begin{aligned}
& \nabla^M_{e_k}\nabla^M_{e_k}(\Omega_2 \wedge \omega)(e_1, \cdots, e_n) \\
& = (\Omega_2 \wedge \nabla^M_{e_k}\nabla^M_{e_k}\omega)(e_1, \cdots, e_n) \\
& = \frac{\lambda_1\lambda_2}{\sqrt{(1+\lambda_1^2)(1+\lambda_2^2)}}(\nabla^M_{e_k}\nabla^M_{e_k}\omega)(a_3, \cdots, a_n)
\end{aligned}$$

For $i = 3, \cdots, n$, $(\nabla^M_X\nabla^M_Y a^*_i)(a_i) = X(\nabla^M_Y a^*_i(a_i)) - (\nabla^M_Y a^*_i)(\nabla^M_X a_i)$. $(\nabla^M_Y a^*_i)(a_i) =< \nabla^M_Y a_i, a_i >=< \nabla^{\Sigma_1}_Y a_i, a_i >$ is zero, so $X(\nabla^M_Y a^*_i(a_i)) = 0$. Therefore $\nabla^M_X\nabla^M_Y a^*_i(a_i) = 0$ at $p$. Since $\omega = a^*_3 \wedge \cdots \wedge a^*_n$, we get $(\nabla^M_{e_k}\nabla^M_{e_k}\omega)(a_3, \cdots, a_n) = 0$.

Now the left hand side of equation (4.10) is non-positive, we shall show the curvature term

$$-\sum_{\alpha,k}[\Omega_{\alpha 2\cdots n} R_{\alpha kk1} + \cdots + \Omega_{1\cdots(n-1)\alpha} R_{\alpha kkn}]$$

is strictly positive, thus achieves contradiction because all other terms on the right hand side are non-negative.

We calculate the curvature term.

$$\Omega(e_\alpha, e_2, \cdots, e_n) = -\delta_{\alpha,n+1}\frac{\lambda_1+\lambda_2}{\sqrt{(1+\lambda^2)(1+\lambda_2^2)}}$$

Likewise,

$$\Omega(e_1, e_\alpha, e_3 \cdots, e_n) = -\delta_{\alpha,n+2}\frac{\lambda_1+\lambda_2}{\sqrt{(1+\lambda^2)(1+\lambda_2^2)}}$$

and

$$\Omega(e_1, e_2, e_\alpha, e_4 \cdots, e_n) = 0$$

We assume $\Sigma_1$ and $\Sigma_2$ are of constant curvature $k_1$ and $k_2$ respectively, by the calculation in [10] we have

$$\begin{aligned}
&\sum_k R(e_\alpha, e_k, e_k, e_i)\\
&= k_1[\sum_k < \pi_1(e_\alpha), \pi_1(e_k) >< \pi_1(e_k), \pi_1(e_1) > - < \pi_1(e_\alpha), \pi_1(e_1) > \sum_k < \pi_1(e_k), \pi_1(e_k) >]\\
&+ k_2[\sum_k < \pi_2(e_\alpha), \pi_2(e_k) >< \pi_2(e_k), \pi_2(e_1) > - < \pi_2(e_\alpha), \pi_2(e_1) > \sum_k < \pi_2(e_k), \pi_2(e_k) >]\\
=&k_1[\sum_k < \pi_1(e_\alpha), \pi_1(e_k) >< \pi_1(e_k), \pi_1(e_i) > - < \pi_1(e_\alpha), \pi_1(e_i) > \sum_k |\pi_1(e_k)|^2]\\
&+ k_2[(n-1) < \pi_1(e_\alpha), \pi_1(e_i) >\\
&+ \sum_k < \pi_1(e_\alpha), \pi_1(e_k) >< \pi_1(e_k), \pi_1(e_1) > - < \pi_1(e_{n+1}), \pi_1(e_1) > \sum_k |\pi_1(e_k)|^2]
\end{aligned}$$

Because $e_{n+i} = \frac{1}{\sqrt{1+\lambda_i^2}}(a_{n+i} - \lambda_i a_i)$, $\pi_1(e_{n+i}) = \frac{-\lambda_1}{\sqrt{1+\lambda_i^2}}a_i$. Likewise, $\pi_1(e_k) = \frac{1}{\sqrt{1+\lambda_k^2}}a_k$.

Therefore,

$$R_{n+i,kki} = \frac{\lambda_i}{1+\lambda_i^2}[k_1(\sum_{k\neq i}\frac{1}{1+\lambda_k^2}) + k_2(1-n+\sum_{k\neq i}\frac{1}{1+\lambda_k^2})]$$

Therefore the curvature term $-\sum_{\alpha,k}[\Omega_{\alpha 2\cdots n}R_{\alpha kk1} + \cdots + \Omega_{1\cdots(n-1)\alpha}R_{\alpha kkn}]$ in equation (3.6) becomes

$$\frac{(\lambda_1+\lambda_2)}{\sqrt{(1+\lambda_1^2)(1+\lambda_2^2)}}\sum_{i=1}^{2}\frac{\lambda_i}{1+\lambda_i^2}[k_1(\sum_{j\neq i}^{n}\frac{1}{1+\lambda_j^2}) + k_2(1-n+\sum_{j\neq i}^{n}\frac{1}{1+\lambda_j^2})]$$

$$k_1(\sum_{j\neq i}\frac{1}{1+\lambda_j^2})+k_2(1-n+\sum_{j\neq i}\frac{1}{1+\lambda_j^2}) = \frac{k_1-k_2}{2}(n-1)+\frac{k_1+k_2}{2}(\sum_{j\neq i}\frac{2}{1+\lambda_j^2}+1-n)$$

Since $\lambda_1+\lambda_2>0$ and $k_1-k_2\geq 0, k_1+k_2\geq 0$, we only need to show

$$\sum_{i=1}^n \frac{\lambda_i}{1+\lambda_i^2}(1-n+\sum_{k\neq i}^n\frac{2}{1+\lambda_k^2})\geq 0$$

This is indeed

$$\frac{\lambda_1}{1+\lambda_1^2}(n-3+\frac{2}{1+\lambda_2^2})+\frac{\lambda_2}{1+\lambda_2^2}(n-3+\frac{2}{1+\lambda_1^2})$$

This can be rewritten as

$$(n-2)\frac{(\lambda_1+\lambda_2)(1+\lambda_1\lambda_2)}{(1+\lambda_1^2)(1+\lambda_2^2)}+\frac{(\lambda_1+\lambda_2)(1-\lambda_1\lambda_2)}{(1+\lambda_1^2)(1+\lambda_2^2)}$$

which is strictly positive under the assumption $|\lambda_1\lambda_2|<1$.

Now we turn to long-time existence and convergence. As in Theorem A, we can show

$$|\lambda_1\lambda_2|\leq 1-\delta, \tag{4.11}$$

Recall the equation satisfied by $*\Omega_1$,

$$\frac{d}{dt}*\Omega_1\geq\Delta*\Omega_1+\delta*\Omega_1|A|^2$$
$$+*\Omega_1\{\frac{k_1-k_2}{2}(n-1)[\sum_i\frac{\lambda_i^2}{1+\lambda_i^2}]+\frac{k_1+k_2}{2}[(n-2)\frac{\lambda_1^2+\lambda_2^2+2\lambda_1^2\lambda_2^2}{(1+\lambda_1^2)(1+\lambda_2^2)}+\frac{\lambda_1^2+\lambda_2^2-2\lambda_1^2\lambda_2^2}{(1+\lambda_1^2)(1+\lambda_2^2)}]\} \tag{4.12}$$

In fact, $\sum_i\frac{\lambda_i^2}{1+\lambda_i^2}=\frac{\lambda_1^2+\lambda_2^2+2\lambda_1^2\lambda_2^2}{(1+\lambda_1^2)(1+\lambda_2^2)}$. Since $*\Omega_1=\frac{1}{\sqrt{(1+\lambda_1^2)(1+\lambda_2^2)}}$, we have $1-*\Omega_1^2=\frac{\lambda_1^2+\lambda_2^2+\lambda_1^2\lambda_2^2}{(1+\lambda_1^2)(1+\lambda_2^2)}$. It is not hard to see there exists a constant $c'$ such that $\lambda_1^2+\lambda_2^2-2\lambda_1^2\lambda_2^2\geq c'(\lambda_1^2+\lambda_2^2+\lambda_1^2\lambda_2^2)$ if $|\lambda_1\lambda_2|\leq 1-\delta$.

Therefore,

$$\frac{d}{dt}*\Omega_1\geq\Delta*\Omega_1+\delta*\Omega_1|A|^2+c*\Omega_1(1-*\Omega_1^2) \tag{4.13}$$

for some constant $c > 0$.

As in the proof of Theorem A in [10], this equation implies long time existence by blowing up argument and White's regularity Theorem [11].

By maximum principle, $\min_{\Sigma_t} *\Omega_1 \to 1$ as $t \to \infty$, then we can use the estimate in the proof of Theorem B in [10] to show $\max_{\Sigma_t} |A|^2 \to 0$ and apply Simon's Theorem. We get smooth convergence in this case. In the limit, $*\Omega_1 = 1$ and thus $\lambda_1 = \lambda_2 = 0$.

□

**Corollary** *If $f$ is any smooth map from $S^n$ to $S^2$ and if at each point, the restriction of $df$ to any two dimensional subspace is area decreasing, then $f$ is homotopic to a constant map along the mean curvature flow.*

Again, we remark the condition is equivalent to the comass of $f^*\Omega_2$ is less than one, where $\Omega_2$ is the area form on $S^2$.

# References

[1] K. Ecker and G. Huisken *Mean curvature evolution of entire graphs*, Ann. of Math. (2) 130 (1989), no. 3, 453–471.

[2] K. Ecker and G. Huisken, *Interior estimates for hypersurfaces moving by mean curvature.*, Invent. Math. 105 (1991), no. 3, 547–569.

[3] H. Federer, *Geometric measure theory.* Die Grundlehren der mathematischen Wissenschaften, Band 153 Springer-Verlag New York Inc., New York 1969 xiv+676 pp.

[4] R. Hamilton *Four-manifolds with positive curvature operator.* J. Differential Geom. 24 (1986), no. 2, 153–179.

[5] R. Harvey and H. B. Lawson, *Calibrated geometries.* Acta Math. 148 (1982), 47–157.

[6] G. Huisken, *Asymptotic behavior for singularities of the mean curvature flow*, J. Differential Geom. **31** (1990), no. 1, 285–299.

[7] L. Simon, *Asymptotics for a class of nonlinear evolution equations, with applications to geometric problems.*, Ann. of Math. (2) 118 (1983), no. 3, 525–571.

[8] M-T. Wang, *Mean Curvature Flow of surfaces in Einstein Four-Manifolds* , J. Differential Geom. **57** (2001), no. 2, 301-338.

[9] M-T. Wang, *Deforming area preserving diffeomorphism of surfaces by mean curvature flow* , Math. Res. Lett. **8** (2001), no. 5-6. 651-661.

[10] M-T. Wang, *Long-time existence and convergence of graphic mean curvature flow in arbitrary codimension* , Invent. math. 148 (2002) 3, 525-543.

[11] B. White, *A local regularity theorem for classical mean curvature flow*, preprint, 2000.